\documentclass[preprint]{article}
\usepackage{graphicx} 
\usepackage[left=2cm, right=2cm, top=2cm, bottom=2cm]{geometry}
\usepackage{algorithm}
\usepackage{algpseudocode}
\usepackage{amsmath,amssymb}
\usepackage{comment}
\usepackage{hyperref}
\usepackage{stmaryrd}
\usepackage{bm}
\usepackage{authblk}

\title{Nonnegative tensor train for the multicomponent Smoluchowski equation}
\author[1,2]{Segey Matveev\thanks{Corresponding author: matseralex@cs.msu.ru}}
\author[2,3]{Ilya Tretyak}
\affil[1]{Lomonosov MSU, Faculty of Computational Mathematics and Cybernetics, Moscow, Russia}
\affil[2]{Marchuk Institute of Numerical Mathematics RAS, Moscow, Russia}
\affil[3]{Moscow Institute of Physics and Technology (National Research University), Dolgoprudnii, Russia}
\date{}

\begin{document}
\maketitle{}
\begin{abstract}
    We propose an efficient implementation of the numerical tensor-train (TT) based algorithm solving the multicomponent coagulation equation preserving the nonnegativeness of solution. Unnatural negative elements in the constructed approximation arise due to the errors of the low-rank decomposition and discretization scheme. In this work, we propose to apply the rank-one corrections in the TT-format proportional to the minimal negative element. Such an element can be found via application of the global optimization methods that can be fully implemented within efficient operations in the tensor train format. We incorporate this trick into the time-integration scheme for the multicomponent coagulation equation and also use it for post-processing of the stationary solution for the problem with the source of particles.       
\end{abstract}


\section{Introduction}

Coalescence (coagulation/aggregation) of interacting particles  plays a fundamental role  in many processes in nature and technology \cite{KrapivskyBook, Leyvraz} such as rain droplet formation \cite{Falkovich2002, aloyan1997transport}, rubber production \cite{de2015acid} and synthesis of nanoparticles \cite{boje2017detailed}. Aggregation of the microparticles is essential for the molecular beam epitaxy \cite{kaganer2016nucleation, sabelfeld2017dislocation} and also important for natural crystal growth \cite{Flory, blatz1945note} and blood clotting \cite{anand2006model, filkova2019quantitative}. On the basic level these processes correspond to a sequence of reactions 
$$[u]  + [v] \xrightarrow[~]{K(u,v)} [u+v]$$
with nonnegative symmetric reaction rates $K(u,v)=K(v,u) \geq 0$. The basic governing equations for such set of reactions read as integro-differential Smoluchowski equations \cite{Smo17}
\begin{equation}\label{eq:smol}
    \frac{\partial n(v,t)}{\partial t} = \frac{1}{2}\int_{0}^{v} K(v-u,u) n(v-u,t) n(u,t) du - n(v,t) \int_{0}^{\infty} K(v,u) n (u,t) du.
\end{equation}
These equations describe the evolution of the particle number density function $n(v,t)$ changing due to \textit{birth} of particles of size $v$ due to the agglomeration on of the smaller particles $v-u$ and $u$ (the first term in the right-hand side) and their \textit{death} after coagulation with other particles (the second term). 

This model (and its discrete formulations) is currently rather well studied \cite{Lushnikov1976, PNAS} and there exists a number of reliable and efficient numerical tools its investigation. These tools can be either deterministic (e.g. finite-difference\cite{matveev_fast_2015} and finite-volume methods\cite{singh2021accurate}) or stochastic dealing with some ensemble of particles and direct modeling of the agglomeration events \cite{eibeck2000efficient, sorokin2015monte}. 

In case of the complex particles consisting of different components (e.g. water, carbon and different acids) and their source the basic equation has to be generalized in multidimensional form \cite{palaniswaamy2007direct} 
\begin{eqnarray}
\notag\frac{\partial n(\overline{v}, t)}{\partial t} = \underbrace{\frac{1}{2}\int_0^{v_1} \ldots \int_0^{v_d} K(\overline{v}-\overline{u}; \overline{u}) n(\overline{v} - \overline{u}, t) n(\overline{u}, t) du_1 \ldots du_d}_{\mathcal{L}_1(\overline{v})} -\\ \label{eq:SmoluchowskiMultidim}
\\ \notag - n(\overline{v},t) \underbrace{\int_0^{\infty} \ldots \int_0^{\infty} K(\overline{u}; \overline{v}) n(\overline{u},t) du_1 \ldots du_d}_{\mathcal{L}_2(\overline{v})} + q(\overline{v}, t),
\end{eqnarray}
where $n(\overline{v},t)$ corresponds to the concentration function of the complex particles consisting of a set $(v_1, v_2, \ldots, v_d)$ of components and function $q(v_1, \ldots, v_d, t)$ represents the sources of different particles. 

Minor correction of the initial problem for the finite particle size domain (meaning the sink of the large particles from the system, e.g. due to the gravity)
\begin{eqnarray}
\notag
\frac{\partial n(\overline{v}, t)}{\partial t} = \frac{1}{2}\int_0^{v_1} \ldots \int_0^{v_d} K(\overline{v}-\overline{u}; \overline{u}) n(\overline{v} - \overline{u}, t) n(\overline{u}, t) du_1 \ldots du_d - \\ \label{eq:SmoluchowskiMultidimSource}
\\ \notag - n(\overline{v},t) \int_0^{V_{\max}} \ldots \int_0^{V_{\max}} K(\overline{u}; \overline{v}) n(\overline{u},t) du_1 \ldots du_d + q(\overline{v}, t),
\end{eqnarray}
often leads to stationary solutions $n(\overline{v})$ corresponding to stable particle size distributions for $t \rightarrow \infty$.

The multicomponent coalescence models are useful for description of miscellaneous processes in atmospheric aerosols consisting from dozens of different chemical compounds \cite{aloyan1997transport, palaniswaamy2007direct, sorokin2015monte}. Due to the higher dimensionality this equation is rather tough for theoretical and numerical investigation. Detailed modeling of the multicomponent coagulation with Monte Carlo methods often meets the unpleasant accuracy limitations and requires utilization of tremendously large numbers of particles. Even though these methods are popular in practice \cite{sorokin2015monte, boje2019hybrid} the precise numerical simulation of the whole particle size distribution function $n(\overline{v},t)$ might require significant effort.

Formulation of the efficient deterministic approach for this problem requires utilization of the advanced methods accounting the low-parametric structure of the studied problem. Application of the tensor train (TT) decomposition \cite{smirnov_fast_2016, matveev_tensor_2016, manzini2021nonnegative} seems to be one of recent successes in this direction. The low-rank decomposition of the solution and coagulation kernel in TT-format allows to obtain a reasonable complexity of the time-integration methods (e. g. Runge-Kutta approach). Hence, they allow to study interesting problems with good precision using usual laptops in modest times. We revisit this approach in Section \ref{sec:TT-method}.

However, tensor-based methods for multicomponent coagulation still have a significant drawback: they do not guarantee the elementwise nonnegativity of the constructed solutions $n(\overline{v}, t)$. One may apply the alternating projections methods \cite{sultonov2023low} or the nonnegative matrix factorization methods \cite{manzini2021nonnegative} to resolve this problem. However, the first approach \cite{sultonov2023low} requires elementwise corrections for the whole data during the iterative process and the latter is applicable only for the two-component coagulation problems.

In Section \ref{sec:NTT} we show that the numerical solution in the tensor train format can be fixed within modest additional complexity and low-rank correction of the solution. This correction can be applied effectively for any exact time-stepping procedure as well as just once for the final numerical approximation of the stationary solution (e.g. for the problem of multicomponent coagulation with source). 

In Section \ref{sec:numerical} we show in practice that our ideas can be incorporated for the problems of arbitrary dimensionality with different coagulation kernels. We present the numerical results for three- and four-dimensional problems and conclude that novel approach has a modest additional computational cost without loss of accuracy in comparison with initial method.

\section{Fast algorithm using Tensor-Train}
\label{sec:TT-method}

In this section we provide a brief introduction to the tensor-train based numerical approach solving the multidimensional Smoluchowski equation. A detailed description and validation of this approach can be found in \cite{matveev_tensor_2016}. 

We consider each component $\overline{v}$ on a finite segment $[0, V_{max}]$, assuming that $V_{max}$ is common for all $(v_1,v_2,\ldots,v_d)$. Then we discretize each segment with $N$ equidistant nodes defining a uniform grid with $N^d$ nodes. For each time moment $t$ the solution is represented by $d$-dimensional grid-function $n(\overline{i}, t)$, where $\overline{i}=(i_1,\ldots,i_d)$ is some point on the grid. Thus, utilization of the classical approximation and time-integration schemes on this grid requires exponential number of memory cells for storage and calculation of the solution. Therefore, we can use the tensor-train representation of $d$-dimensional solution at each time moment
\begin{equation}\label{eq:Solution TT-format}
n(v_1,\ldots,v_d,t) = \sum\limits_{\alpha_0=1}^{r_0}\sum\limits_{\alpha_1=1}^{r_1}\ldots\sum\limits_{\alpha_d=1}^{r_d}n_1(\alpha_0, v_1, \alpha_1, t)\ldots n_d(\alpha_{d-1}, v_d, \alpha_d, t),
\end{equation}
where $r_k$ are called ranks and tensor-train representation (TT-format) assumes that $r_0=r_d=1$, $n_\mu$ are the time-dependent tensor-train cores. The major advantage of the TT-format is that we need to store $O(dNr^2)$ memory cells at the each time moment (where $r=\max{r_k}$) instead of original $O(N^{d})$ for the direct method. We can evaluate the concentration function in TT-format at any grid-point in $O(dr^2)$ operations using the sequence of matrix by vector product operations. Summarizing with the fact that coagulation kernels and solution of the Smoluchowski equation have small ranks for a wide class of functions \cite{matveev_tensor_2016}, we can build a very memory-efficient algorithm having modest algorithmic complexity.

We use the classical predictor-corrector scheme for the time integration. Fixing the segment $[0, T]$, we introduce a virtual time-grid with step $\tau$ and denote $n^k(\overline{i})=n(\overline{i}, k\tau)$. Thus, the time integration scheme has the form
\begin{eqnarray}
\notag
    n^{k+\frac{1}{2}}\left(\overline{i}\right) = \frac{\tau}{2} \left(L_1^k\left(\overline{i}\right) - n^k\left(\overline{i}\right)L_2^k\left(\overline{i})\right)\right) + n^k\left(\overline{i}\right),\\
\label{eq:Discrete approxiamtion}
    \\ \notag n^{k+1}\left(\overline{i}\right) = \tau\left(L_1^{k+\frac{1}{2}}\left(\overline{i}\right) - n^{k + \frac{1}{2}}\left(\bar{i}\right)L_2^{k+\frac{1}{2}}\left(\bar{i}\right)\right) + n^k\left(\bar{i}\right),
\end{eqnarray}
where upper index in integral-terms approximations denotes the solution at the corresponding time-step. Note that application of this scheme is simply obtained by using only the TT-arithmetic operations because all data is stored in this format. Thus, we can work exclusively with tensor-train representation throughout the integration time.

We approximate the integral operators $\mathcal{L}_1$, $\mathcal{L}_2$  from the original equation \eqref{eq:SmoluchowskiMultidim} by multidimensional trapezoid rule and denote the corresponding grid versions $L_1(\overline{i})$ for the first integral and $L_2(\overline{i})$ for second. Direct computation of this approximation leads to $O(N^{2d})$ operations requiring significant computing resources. 

Hence, we use modification of the fast algorithm for calculation of the lower-triangular convolution in tensor-train format for $L_1(\overline{i})$ which costs $O(d^2R^4N\log N)$ operations, where $R$ is maximum rank of convolution members. Reminding the original idea from \cite{matveev_tensor_2016} we write them as following
\begin{eqnarray}
    \notag \mathcal{L}_1 (\overline{v}) = \sum\limits_{\alpha_0,\ldots,\alpha_{2d}}\sum\limits_{\beta_0\ldots,\beta_d} \int\limits_{0}^{v_1}{K_1^v(\alpha_0,v_1-u_1,\alpha_1)K_1^u(\alpha_d, u_1, \alpha_{d+1})n_1(\beta_0,v_1-u_1,\beta_1,t)n_1(\beta_0,u_1,\beta_1,t)du_1}\times \ldots\\
    \label{eq:first_tt_integral}
    \\ \notag \ldots\times \int\limits_{0}^{v_1}{K_d^v(\alpha_{d-1},v_d-u_d,\alpha_d)K_d^u(\alpha_{2d-1}, u_d, \alpha_{2d})n_d(\beta_{d-1},v_d-u_d,\beta_d,t)n_d(\beta_{d-1},u_d,\beta_d,t)du_d}.
\end{eqnarray}
Each integral along $u_i$ mode can be computed as classical convolution via the Fast Fourier Transform (FFT), elementwise product and inverse FFT operations \cite{matveev_tensor_2016}. Further, we re-write the second integral operator
\begin{eqnarray}
    \notag\mathcal{L}_2 (\overline{v})= \sum\limits_{\alpha_0,\ldots,\alpha_{2d}}K_1^v(\alpha_0,v_1,\alpha_1)\ldots K_d^v(\alpha_{d-1},v_d,\alpha_{d}) \sum\limits_{\beta_0\ldots\beta_d} \int\limits_0^{V_{max}}K_1^u(\alpha_d, u_1, \alpha_{d+1})n_1(\beta_0,u_1,\beta_1,t)du_1 \times \ldots \\
    \label{eq:second_tt_integral}
    \\ \notag \ldots \times \int\limits_0^{V_{max}}K_d^u(\alpha_{2d-1}, u_d, \alpha_{2d})n_d(\beta_{d-1},u_d,\beta_d,t)du_d.
\end{eqnarray}
Its grid version  $L_2(\overline{i})$ may be evaluated with second order accuracy taking $O(d^2NR^4)$ operarions using standard TT-arithmetic operations by trapezoid rule . It is worth saying that in practice, there is no need to afraid of a quadratic dependency on the dimension of $d$, because it arises from the application of the TT-cross method and the computed TT-ranks decrease in comparison with the estimate. All in all, we get a semi-discrete approximation, which right-hand side can be evaluated without need to leave the compressed TT-format
\begin{equation}\label{eq:Semi-discrete approxiamtion}
    \frac{\partial n(\overline{i},t)}{\partial t} = \frac{1}{2}L_1(\overline{i}) - n(\overline{i},t) \cdot L_2(\overline{i}).
\end{equation}
The lack of nonnegativeness of the solution due to the inevitable accumulation of the approximation errors during the computations is a certain drawback of this methodology. In the next section we propose an algorithm solving this problem. 

\section{Nonnegative correction}
\label{sec:NTT}

Errors of the tensor-train rank truncation may be accumulated within a few time steps and lead to negative concentration values, which are very unpleasant for the modelling purposes. In this section, we offer a way to deal with occurrence of negative values. Our approach is based on idea nonnegative tensor-train corrections which was introduced in \cite{shcherbakova2022fast, shcherbakova2023study}. This method allows to get a correction of all negative elements for a tensor train of arbitrary dimensionality in contrast with earlier HALS-based method for the two-dimensional case \cite{manzini2021nonnegative}. Our method also allows to do all computations in the compressed TT-format in contrast to other recent approaches \cite{matveev2023sketching, sultonov2023low, jiang2023nonnegative} utilizing the alternating projections.

Let $n^k$ be the solution at $k$-th time step in tensor train format. If the corresponding tensor has no negative values then there is no need in the correction. However, if the negative values exist then we correct the solution by adding the rank-one tensor with elements equal to the absolute of minimal value 
\begin{equation}\label{eq:Nonegative correction}
    \widetilde{n}^k = n^k + |\min_{\overline{i}} n^k(\overline{i})|\cdot \mathbf{E},
\end{equation}
where $\mathbf{E}$ is a tensor with all elements equal to one (some theory about the rank-one corrections can be found e.g. in  \cite{budzinskiy2023quasioptimal}). The result is a tensor-train corresponding to the nonnegative grid-function. Note that, we do not apply a sequence of pointwise rank-one sparse corrections as it has been proposed in \cite{shcherbakova2022fast, shcherbakova2023study}. We do the only rank-one shift because a lot of negative elements in our data are close to the absolute minimum. 

There are various ways allowing to estimate the minimal element of tensor in the TT-format. In our work we use the popular implementation of the global TT-optimization \cite{zheltkov2020global, sozykin2022ttopt} from the TT-toolbox package in Matlab \cite{oseledets2016tt}. One can also apply the iterative approaches, e.g. the power iteration utilizing various acceleration tricks \cite{lebedeva2011tensor}. We utilize the procedure of global optimization in TT-format having the complexity $O(dNR^3)$ operations. Algorithmic cost of this operation does not exceed the total complexity of time-integration step for the multicomponent coagulation equation. Hence, application of the nonnegative correction does not lead to a heavy increase computation time. 

As a result of global optimization we can evaluate absolute maximum of the solution. In order to estimate the minimum, we apply this procedure for the shifted tensor 
$$n^k - M\cdot \mathbf{E},$$
where $M$ is estimated maximum. Thus, we get minimum approximation
\begin{equation}\label{eq:min}
    |\min_{\overline{i}} n^k(\overline{i})|\approx |m - M|,
\end{equation}
where $m$ is maximum found via application of TT-optimization for a second time.

\begin{algorithm}[ht!]
\caption{Nonnegative tensor-train based solver}\label{solver}
\begin{algorithmic}[1]
\Procedure{Nonnegative solver}{}

Get TT-representations of initial condition $n^0(\overline{i})$ and coagulation kernel $K$,

$\tau$ - time step,

$N_t$ - number of time steps,

$s$ - nonegative correction step

\While{$k \leq N_t$}
\State $n^{k}$ $\gets$ $do\_predictor\_corrector\_step$ $(n^{k-1}, K, \tau)$

\If{$k == s$}
\State $M \gets$ $TT\_max(n^k)$ \Comment{search maximum}
\State $m \gets$ $TT\_max(n^k - M\cdot I)$ \Comment{search abs(minimum - maximum)}
\If{$(m - M) > 0$} \Comment{checking for negative elements}
    \State $n^k = n^k + (m - M)\cdot TT\_ones(n^k)$ \Comment{do correction}
\EndIf
\EndIf

\EndWhile

\EndProcedure
\end{algorithmic}
\end{algorithm}

Such a correction is supposed to be applied after the time steps when negative elements appear in the approximate solution or once for post-processing of the final data in calculations (e.g. for the multicomponent coagulation with a source). In practice, checking presence of negative elements in a multidimensional solution at each time step would imply their search, which constitutes the main part of computational complexity for nonnegative correction procedure. For this reason, we correct the numerical solution after a certain pre-selected number of steps. 

We present the final pseudocode of our approach with Algorithm \ref{solver}. The main advantage of Algorithm \ref{solver} is that the solution is stored within the tensor train format throughout full integration time allowing us to work with high-dimensional data.

\section{Numerical experiments}
\label{sec:numerical}

In this section we present our numerical experiments for the nonnegative approximation of the Smoluchowski equation and discuss the modelling in more detail. We compare the proposed nonnegative approximation methodology with older approach using tensor-train selecting it as a benchmark for the fast solver. We study the usual multicomponent coagulation equation \eqref{eq:SmoluchowskiMultidim} and source-enhanced problem \eqref{eq:SmoluchowskiMultidimSource}. In all cases we setup exponential inital conditions 
\begin{equation*}
    n(v_1,v_2, \ldots, v_d,t=0) = e^{-v_1-\ldots-v_d}
\end{equation*}
for the corresponding Cauchy problems.

We present our results for the two different coagulation kernels:
\begin{itemize}
    \item The simplest possible size-independent kernel 
\begin{equation}
    \label{eq:const_kernel}
    K(\overline{u}, \overline{v}) = 1
\end{equation}
with a known analytical solution for the exponential initial conditions and $d=2$. This kernel has a trivial TT-decomposition with all TT-ranks equal to one.
\item Complicated ballistic kernel without known exact solution 
\begin{equation}
    \label{eq:ballistic_kernel}
    K(\overline{u}, \overline{v}) = \left (\left( \sum\limits_{i=1}^d{v_i} \right)^{1/3} + \left( \sum\limits_{i=1}^d{v_i} \right)^{1/3} \right)^2 \sqrt{\frac{1}{\sum\limits_{i=1}^d{u_i}} + \frac{1}{\sum\limits_{i=1}^d{v_i}}}
\end{equation}
having more complex but still low-rank TT-representation (see e.g. \cite{matveev_tensor_2016}).
\end{itemize}
For all numerical experiments we utilized the laptop with Intel(R) Core(TM) i5-8257U CPU @ 1.40GHz processor with 8 Gb RAM and the matlab implementation of TT-toolbox \cite{oseledets2016tt}.

\begin{figure}[ht!]
\begin{center}
        \includegraphics[scale=0.1]{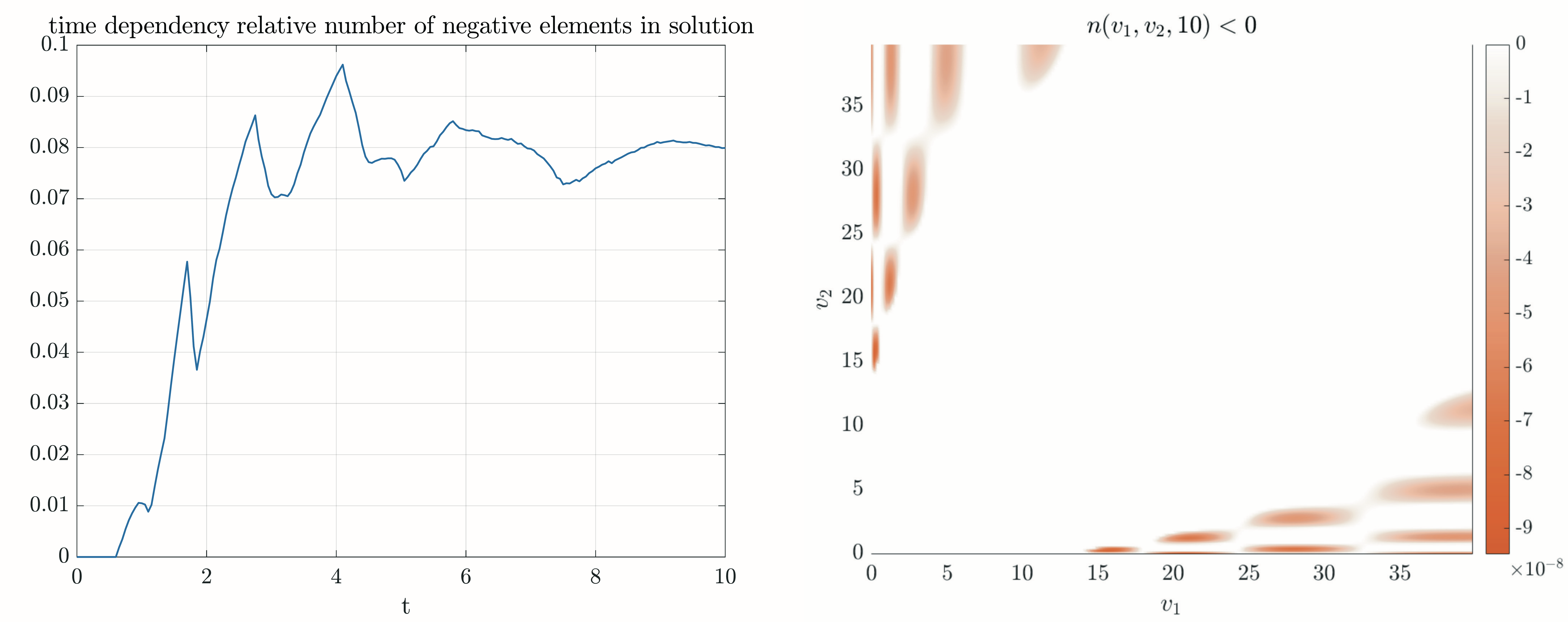}
        \caption{The two-component problem \eqref{eq:SmoluchowskiMultidim} with the constant kernel and exponential initial conditions with $h=0.1, \tau=0.05$. Left panel: Dynamics of the relative number of negative elements; Right panel: negative elements in the numerical solution at $t=10$. }
        \label{fig:negative_elements}
\end{center}
\end{figure}

\subsection{Coagulation without source}

At first, we study the multicomponent coagulation equation \eqref{eq:SmoluchowskiMultidim} without source. The case of the constant kernel with initial conditions
\begin{equation}
    n(v_1,v_2,0) = ab \ e^{-av_1-bv_2}, \ a,b>0
\end{equation}
is exactly solvable and its analytical solution is well-known (see e.g. \cite{fernandez2007exact}) and reads
\begin{equation}
    \label{eq:analitical}
    n(v_1,v_2,t) = \frac{ab \ e^{-av_1-bv_2}}{(1+t/2)^2} I_0\left( 2\sqrt{\frac{ab \ v_1v_2t}{t+2}} \right), 
\end{equation}
where $I_0$ is the modified Bessel function of zeroth order. Fixing $a = b = 1$ for this solution one may compute the total density of the particles per the unit volume
\begin{equation}
    \label{eq:analitical_total_density}
    N(t) = \int\limits_0^\infty\int\limits_0^\infty n(v_1,v_2,t)dv_1dv_2 = \frac{2}{2 + t}.
\end{equation}
Mass conservation law is natural for a basic coagulation model \eqref{eq:SmoluchowskiMultidim} without source
\begin{equation}
M(t) = \int_{0}^{\infty} \ldots \int_{0}^{\infty} (v_1 + v_2 + \ldots v_d) ~ n(v_1, \ldots, v_d,t)~ dv_1 \ldots dv_d \equiv \text{const}
\end{equation}
and takes place for a broad family of homogeneous coagulation kernels $K(x \cdot \overline{u},x \cdot \overline{v}) = x^{\nu} \cdot K(\overline{u},\overline{v})$ with $\nu \leq 1$ \cite{melzak1957scalar, melzak1957scalarPt2}. 

In case of of the model \eqref{eq:SmoluchowskiMultidimSource} with sources and sinks of the particles the total mass $M(t)$ does not stay constant. Its dynamics is defined by the balance between the injection and loss of the matter in the system. However, tracking of the total mass $M(t)$ might be useful for analysis of the relaxation of the whole solution $n(\overline{v},t)$ to the stationary state as we demonstrate below in Figure \ref{fig:total_characteristics_source}.

\begin{figure}[ht!]
\begin{center}
        \includegraphics[scale=0.1]{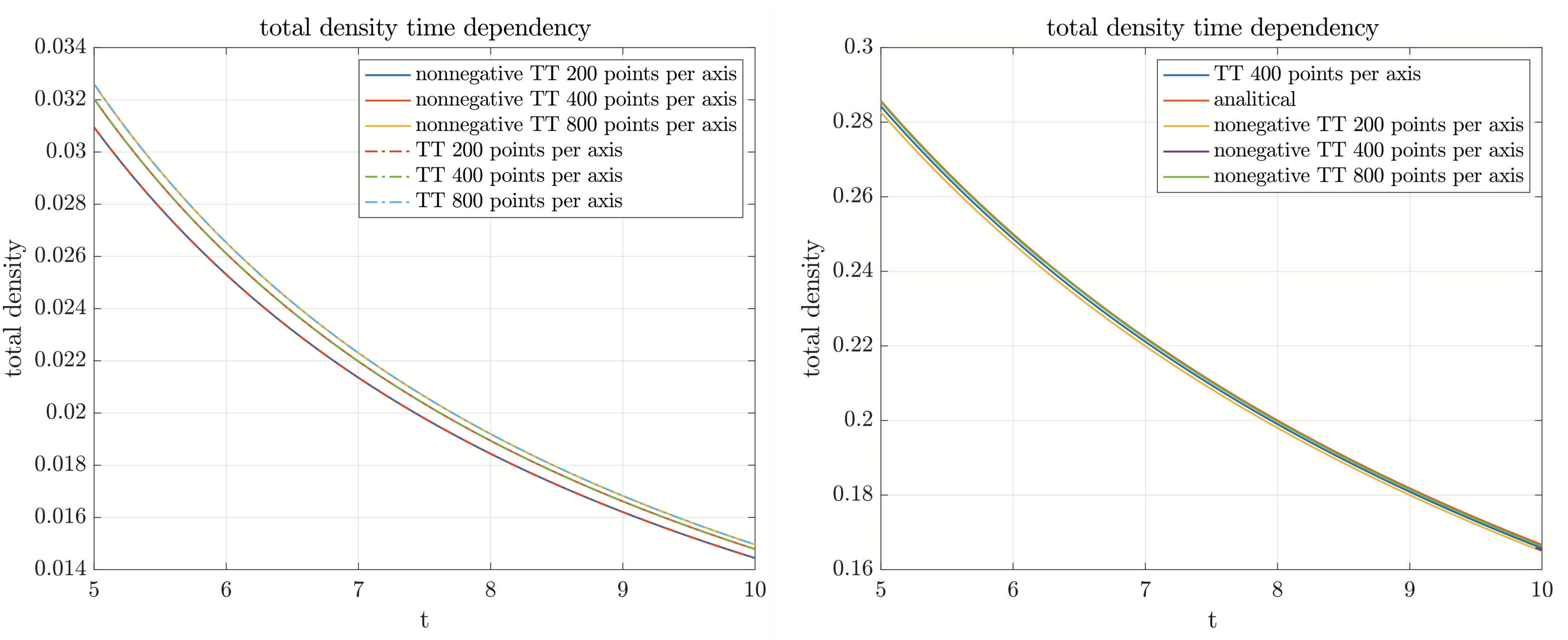}
    \caption{Plots of the total density for the two-component problem with the ballistic kernel (left) and constant kernel (right) and exponential initial conditions for $t\in[5,10]$ simulations time with $V_{max} = 20$. The major difference between the numerical total density $N(t)$ corresponds to the convergence of the finite-difference scheme but not to the introduction of the nonnegative corrections.}
    \label{fig:total_density}
\end{center}
\end{figure}

In Figure \ref{fig:negative_elements} we show the numerical results for the standard TT-based solution without nonnegative corrections. The most significant positive concentrations in the numerical solution $n(v_1, v_2,t)$ have the magnitude $10^{-3}$. On the left panel of we show the dynamics of the relative number of negative elements in the numerical solution. The right panel of corresponds to their heatmap for $t=10$.  One can see that the absolute value of negative elements is quite small and corresponds to the order $10^{-7}$. However, they are clearly visible and their relative number grows up to almost $10 \%$. Hence, we seek to avoid accumulation of such elements without loss of the TT-based compression of the solution. 

In Figure \ref{fig:total_density} we compare the dynamics of the total density for nonnegative TT-approach, basic TT-approach and analytical solution in case of the constant kernel. For each experiment the nonnegative correction is applied after equal fraction of time steps as soon as negative elements arise in the numerical solution. The resulting data generated by the standard approach and by the novel one are almost same in terms of accuracy for the total density.

\begin{table}[ht!]
\begin{center}
\caption{Convergence for the nonnegative TT-based method and the standard TT-based for the 2-dimensional problem with $V_{max} = 40$ and $t \in [0,5]$}.
\label{tab:convergence}
\begin{tabular}
{p{0.5cm}|p{0.75cm}|p{0.75cm}|p{3cm}|p{3cm}|p{3cm}|p{3cm}} 
\hline
N & $\tau$ & $\varepsilon_{cross}$ & Relative error in $\|\cdot\|_F$ for nonnegative TT-solution & Relative error in $\|\cdot\|_F$ for standard TT-solution & Relative error of total density for nonnegative TT-solution & Relative error of total density for standard TT-solution\\
\hline
 200  & 0.1    &1e-6& 5.8605e-2  & 5.8473e-2  & 2.012e-2  & 2.0365e-2\\ 
 400  & 0.05   &1e-7& 2.8363e-2  & 2.8303e-2  & 9.5281e-3 &  9.6054e-3\\ 
 800  & 0.025  &1e-8& 1.3925e-2  & 1.3921e-2  & 4.6246e-3 & 4.6396e-3\\
 1600 & 0.0125 &1e-9& 6.897e-03  & 6.8945e-3  & 2.2868e-3 & 2.2908e-3\\
\hline
\end{tabular}
\end{center}
\end{table}

In Table \ref{tab:convergence} we summarise the convergence results to the exact solution by the nonnegative TT-based solver and compare its accuracy with the standard TT-solver. Periodic nonnegative corrections keep cleanup of the solution from the artefacts and have a very modest influence on the final accuracy of the whole numerical scheme in terms of relative Frobenius norm and relative total density of concentration.

Surprisingly, in case of the modest dimensionality of the corresponding dense tensor corresponding to the solution it \textit{can be} directly placed in memory. In this case, the minimal element can be evaluated accurately with least cost using standard functions of well-optimized linear algebra packages. Thus, a straight-forward search of such an element within $N^{d}$ operations might be faster than application of many TT-power iterations for relatively small values of $N$ and $d$. 

Such a trick allows to abstract from minimum searching speed of particular method and concentrate only on performance of our approach. All the rest of operations (multidimensional integral operators and time-integration steps) must be done with use of the compressed format -- complexity of their brute force evaluation grows as $N^{2d}$ and makes the naive computations impractical even for $N=200$, $d=2$ (see e.g. \cite{matveev_tensor_2016}). Hence, utilization of tensor trains for numerical integration of the Smoluchowski equation requires much less time than the direct method even for the small data. 

In Table \ref{tab:const_table} we show the computational times for the standard TT-based methodology and the novel methodology in exact minimum search for $d=2,3$ and global optimization approach from \cite{zheltkov2020global, sozykin2022ttopt} for $d=4$. In Table \ref{tab:ballistic_table} present the similar measurements for the ballistic kernel. The simulation time increases significantly for the ballistic kernel due to the higher ranks of its TT-representation. 

In addition to evaluation of the relative accuracy of the novel approach we verify its robustness with respect to different values of the convergence criteria $\varepsilon_{cross}$ of the TT-cross approximation method (see Tables \ref{tab:convergence} -- \ref{tab:speed_sorce_table}). We vary $\varepsilon_{cross}$ from $10^{-6}$ to even $10^{-10}$ setting the internal approximations error by the TT-cross operation less than error of original the finite-difference method. The maximal values of the TT-ranks of the solutions  also stay small during all simulations utilizing either the baseline or the new nonnegative method.

\begin{table}
\begin{center}
\caption{Performance of the nonnegative TT-based method and standard TT-based approach for the coagulation equation with the constant kernel, $V_{max} = 40$ and $t\in[0,1]$.}
\label{tab:const_table}
    \begin{tabular}{p{0.75cm}|p{0.75cm}|p{1cm}|p{0.9cm}|p{2.5cm}|p{2cm}|p{2.5cm}|p{2cm}|p{2.2cm}}
         \hline 
         $d$ & $N$ & $\tau$ & $\varepsilon_{cross}$ & Nonneg. TT-method, sec & max NTT-rank & TT-method, sec & max TT-rank & $\dfrac{\|NTT - TT\|_F}{\|TT\|_F}$\\
         \hline
         2 & 200 & 0.1      & 1e-6 &  1.4  &  6  & 1.5   & 5 & 6.8e-8\\
         2 & 400 & 0.05     & 1e-7 & 4.1   &  7  & 3.4   & 6 & 9.1e-8\\
         2 & 800 & 0.025    & 1e-8 & 12.6  &  7  &  11.9 & 6 & 1.1e-6\\
         2 & 1600 & 0.0125  & 1e-9 & 43.6  &  8  & 41.7  & 7 & 1.1e-8\\
         2 & 3200 & 0.00625 & 1e-10 & 170.1 &  9  & 163.8 & 8 & 9.8e-10\\
         \hline
         3 & 200 & 0.1      & 1e-6  & 4.7    & 6 & 4.6    & 5 & 3.1e-5\\
         3 & 400 & 0.05     & 1e-7  &  20.2  & 7 & 19.7   & 6 & 8.5e-6\\
         3 & 800 & 0.025    & 1e-8  &  112.7 & 8 & 82.9   & 7 & 5.5e-6\\
         3 & 1600 & 0.0125  & 1e-9  & 474.5  & 8 & 362    & 7 & 1.7e-7\\
         3 & 3200 & 0.00625 & 1e-10 & 1823.6 & 9 & 1807.1 & 8 & 8.6e-9\\ 
         \hline
         4 & 200 & 0.1       & 1e-6  & 10.2  & 6 & 9       & 5 & 5.9e-4\\
         4 & 400 & 0.05      & 1e-7  & 41.5  & 7 & 33.2    & 6 & 1.1e-4\\
         4 & 800 & 0.025     & 1e-8  & 290.5 & 8 & 199.6   & 7 & 3.5e-5\\
         4 & 1600 & 0.0125   & 1e-9  & 859   & 9 & 849.5   & 8 & 9.9e-6\\
         4 & 3200 & 0.00625  & 1e-10 & 4823  & 9 & 4591.7  & 8 & 6.5e-7 \\
         \hline
    \end{tabular}
\end{center}
\end{table}

\begin{table}
\begin{center}
\caption{Comparison of the nonnegative TT-based method and standard TT-based scheme for the coagulation equation with the ballistic kernel, $V_{max} = 40$ and $t\in[0,1]$.}
\label{tab:ballistic_table}
    \begin{tabular}{p{0.75cm}|p{0.75cm}|p{0.75cm}|p{0.75cm}|p{2.5cm}|p{2cm}|p{2.5cm}|p{2cm}|p{2.2cm}}
         \hline 
         $d$ & $N$ & $\tau$ & $\varepsilon_{cross}$ & Nonneg. TT-method, sec & max NTT-rank & TT-method, sec & max TT-rank & $\dfrac{\|NTT - TT\|_F}{\|TT\|_F}$ \\
         \hline
         2 & 200 & 0.1   &1e-6& 70.5   & 13 & 70.3  & 12 & 5.1e-4\\
         2 & 400 & 0.05  &1e-7& 381.9  & 16 & 370.6 & 15 & 5e-5\\
         2 & 800 & 0.025 &1e-8& 1736.5 & 19 & 1653  & 18 & 8.8e-6\\
         \hline
         3 & 100 & 0.2   &1e-6& 363.1  & 12 & 204.4  & 11 & 8.1e-4\\
         3 & 200 & 0.1   &1e-7& 1008.3 & 14 & 971.4  & 13 & 8.9e-5\\
         3 & 400 & 0.05  &1e-8& 8549.9 & 20 & 7916.1 & 18 & 2.9e-5\\
         \hline
         4 & 100 & 0.2   &1e-6& 422.9   & 15 & 388.9   & 14 & 2.1e-3 \\
         4 & 200 & 0.1   &1e-7& 5184.9  & 25 & 4945.3  & 23 & 6.2e-4\\
         4 & 400 & 0.05  &1e-8& 26797.6 & 49 & 23562.3 & 46 & 2.8e-4\\
         \hline
    \end{tabular}
\end{center}
\end{table}

\subsection{Source-induced coagulation}

We also present the experiments for the Smoluchowski equation \eqref{eq:SmoluchowskiMultidimSource} with exponential source
\begin{equation}
    \label{eq:source}
    q(v_1,\ldots, v_d,t) = e^{-v_1-\ldots-v_d}
\end{equation}
and the same exponential initial conditions. In this case the solution relaxes to the stationary mode due to the balance between the source of particles and their sink through the boundary $V_{max}$ \cite{smirnov_fast_2016}. 

\begin{figure}[ht!]
\begin{center}
        \includegraphics[scale=0.1]{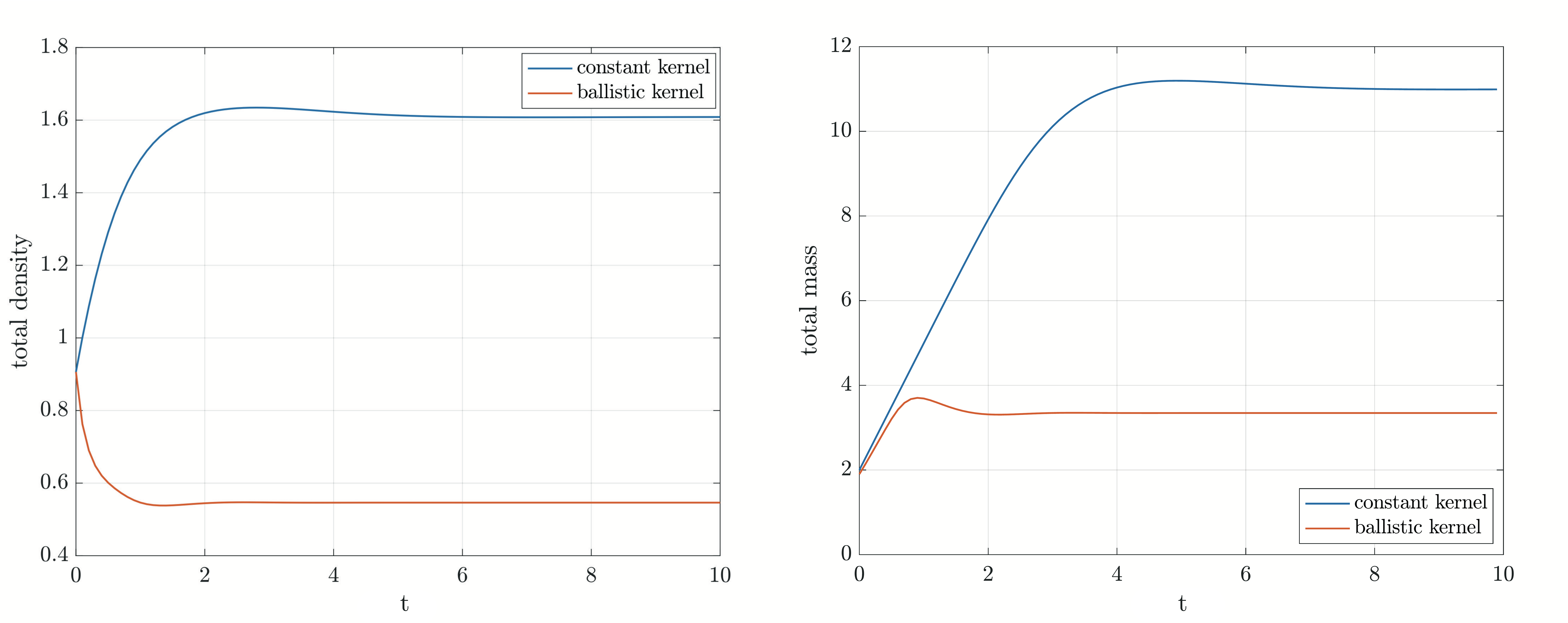}
        \caption{Total density (left) and mass (right) for experiments with source for constant and ballistic kernels for $V_{max}=20$ and $t\in[0,10]$.}\label{fig:total_characteristics_source}
\end{center}
\end{figure}

In Figure \ref{fig:total_characteristics_source} we present the dynamics of the total mass and total density for the constant and ballistic coagulation kernels. We approximate the stationary solution stopping the dynamic computations for some relatively large finite moment (in our simulations we set $T=10$). Hence, it is enough to make only the nonnegative correction just once at the final integration time or after entrance to stationary mode. Such a trick saves computational resources and requires less memory. 

In Table \ref{tab:speed_sorce_table} we present our benchmarks and the relative difference between the solutions constructed by the nonnegative and baseline approaches in the Frobenius norm ($||  \cdot||_F$) for the ballistic coagulation kernel. We demonstrate the results for both implementations of the TT-approach: with nonnegative correction and without. The integration time for the nonnegative case naturally increases due to the efforts for the evaluation of the minimal element in the TT-representation of solution. We see that application of the nonnegative correction requires only very modest additional computational time in comparison with cost of the whole time-integration procedure.  

However, the increase of the TT-ranks of the solution leads to the significant increase of simulation time in agreement with theoretical complexity estimates (see e.g. Table \ref{tab:ballistic_table}). Luckily the application of the nonnegative corrections has almost no impact on the values of the TT-ranks of the solutions.

\begin{figure}[ht!]
\begin{center}
        \includegraphics[scale=0.1]{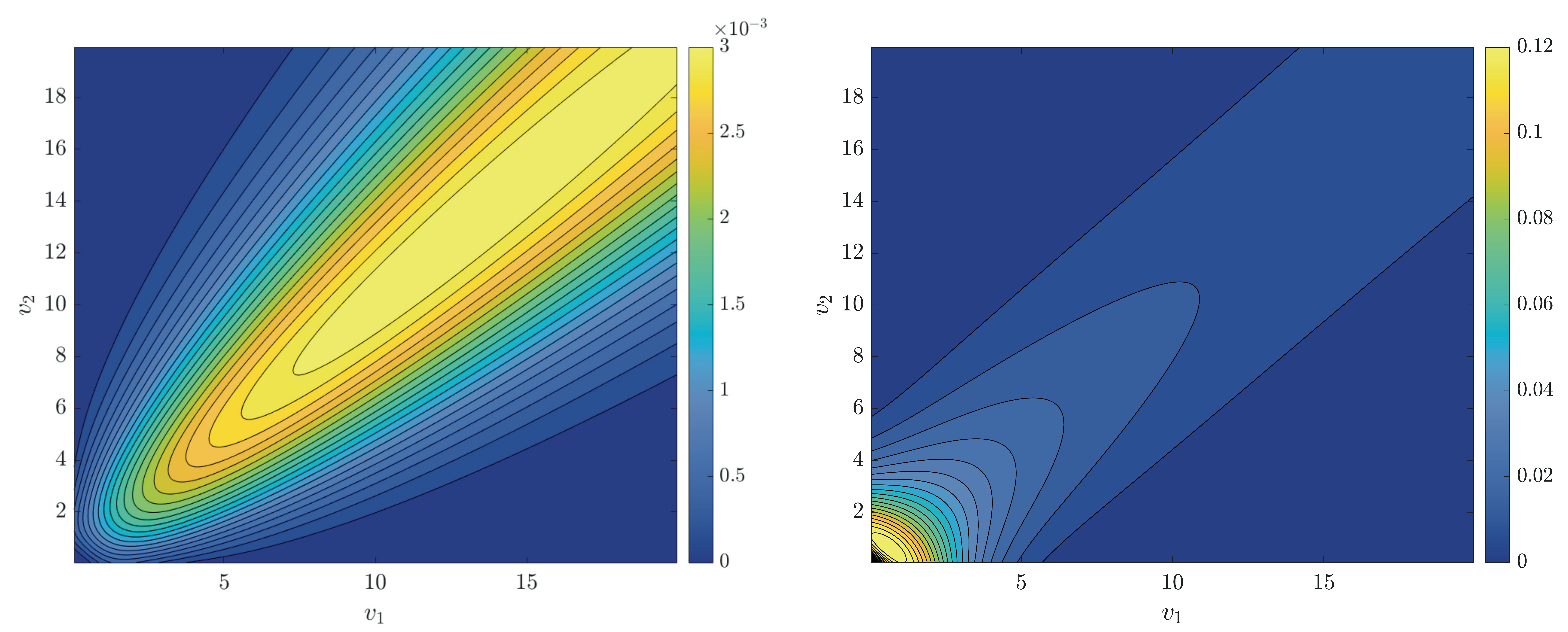}
        \caption{Nonnegative solution for the two-component coagulation problem with the ballistic kernel  without (left) and with source (right), $V_{max}=20$ and $t=5$.}\label{fig:solution_ballistic}
\end{center}
\end{figure}

In fact, the reason of the increase of the ranks is very natural. One may see the plots of the numerical solutions in case of the ballistic kernel in Fig. \ref{fig:solution_ballistic} for $d = 2$. We see that most of \textit{noticeable} elements concentrate and stretch along the diagonal of the computational domain with fast decrease of their amplitudes further the relatively narrow band. This fact allows us to assume that the solution has good structure for its possible compression in the mosaic-skeleton format. However, application of the mosaic skeleton approximations cannot be utilized in the straight-forward way: it would be necessary to reformulate the effective algorithms for computation of the coagulation integrals. We think that such a study even for $d=2$ might lead to fruitful results in future.

\begin{table}[ht!]
\begin{center}
\caption{Comparison for nonnegative and standard TT-methodologies with exponential source for the ballistic kernel, $V_{max} = 20$ and $t\in[0,5]$.}
\label{tab:speed_sorce_table}
    \begin{tabular}{p{0.75cm}|p{0.75cm}|p{0.75cm}|p{0.75cm}|p{4cm}|p{3cm}|p{3cm}}
         \hline 
         $d$ & $N$ & $\tau$ & $\varepsilon_{cross}$ & Nonneg. TT-method, sec & TT-method, sec & $\dfrac{\|NTT - TT\|_F}{\|TT\|_F}$ \\        
         \hline
         2 & 100 & 0.2   &1e-6& 136& 135& 3.8e-6\\
         2 & 200 & 0.1   &1e-7& 542.3& 542.2& 1e-6\\
         2 & 400 & 0.05  &1e-8& 3 539.4& 3 539.4& 1.2e-7\\
         \hline
         3 & 100 & 0.2   &1e-6& 1 463.3& 1 463.2& 4.1e-5\\
         3 & 200 & 0.1   &1e-7& 9 239.2& 9 239.1& 8.5e-6\\
         3 & 400 & 0.05  &1e-8& 164 058.9& 164 058.8& 4.1e-6\\
         \hline
    \end{tabular}
\end{center}
\end{table}

\section{Conclusions}

In this work we demonstrate an application of the nonnegative low-rank tensor train decomposition to numerical solution of multicomponent coagulation equation. Our idea allows to exploit the structure of the compressed solution and construct the nonnegative correction without processing of all elements of the decomposed solution. Such a trick allows us to consider the problems of dimensionality $d > 2$. We test our approach and verify its accuracy for the problems with constant kernel and known analytical solutions. At the same time, we demonstrate that our ideas are applicable to problems with more complicated ballistic kernel and source-induced problems. 

There are still a lot of possible directions for its development: performance of the TT-max procedure might be studied in more detail with use of novel heuristics \cite{sozykin2022ttopt} as well as one may test more advanced methods for finding the nonnegative tensor train decomposition \cite{shcherbakova2022fast, shcherbakova2023study}.  Another interesting challenge is utilization of the virtual dimensions allowing to represent the solution in QTT format \cite{timokhin2020tensorisation} in combination with nonnegative factorization methods. All in all, the coagulation model itself can be generalized by accounting the spontaneous fragmentation terms as well as the coefficients for the charged particles \cite{smith1999coagulation}. We hope to investigate these problems in our future research.

Sergey Matveev was supported by the Russian Science Foundation (project no.  \href{https://www.rscf.ru/project/21-71-10072/}{21-71-10072}). Ilya  Tretyak was supported by Moscow Center of Fundamental and Applied Mathematics at INM RAS, Ministry of Education and Science of the Russian Federation (Grant No. 075-15-2022-286).

\bibliography{NTTSmol}
\bibliographystyle{ieeetr}

\end{document}